\documentclass[preprints,article,accept,moreauthors,pdftex]{mdpi} 

\setitemize{parsep=6pt,itemsep=0pt,leftmargin=*,labelsep=5.5mm}
\setenumerate{parsep=6pt,itemsep=0pt,leftmargin=*,labelsep=5.5mm}
\setlist[description]{itemsep=0mm} 

\firstpage{1} 
\pubvolume{xx}
\issuenum{1}
\articlenumber{5}
\pubyear{2020}
\copyrightyear{2020}
\history{Received: 28 February 2020; Accepted: 27 April 2020; Published: date}
\updates{yes} 

\Title{On the Non 
Metrizability of Berwald 
Finsler~Spacetimes}

\Author{Andrea Fuster $^{1,\dagger}$, Sjors Heefer $^{1,\dagger}$, Christian Pfeifer $^{2,}{}$*$^{,\dagger}{}$ and Nicoleta Voicu $^{3,\dagger}$}

\AuthorNames{Andrea Fuster, Sjors Heefer, Christian Pfeifer and Nicoleta Voicu}

\address{%
$^{1}$ \quad Department of Mathematics and Computer Science, Eindhoven University of Technology, Groene Loper 5, 
5612AZ Eindhoven, The~Netherlands; a.fuster@tue.nl (A.F.); s.j.heefer@tue.nl (S.H.)\\

$^{2}$ \quad Laboratory of Theoretical Physics, Institute of Physics, University of Tartu, W. Ostwaldi 1, 50411 Tartu, Estonia; christian.pfeifer@ut.ee\\
$^{3}$ \quad Faculty of Mathematics and Computer Science, Transilvania University, Iuliu Maniu Str. 50, 500091 Brasov, Romania; nico.voicu@unitbv.ro
}

\corres{Correspondence: christian.pfeifer@ut.ee}

\firstnote{These authors contributed equally to this work.}

\abstract{We investigate whether Szabo's metrizability theorem can be extended to Finsler spaces of indefinite signature. For smooth, positive definite Finsler metrics, this important theorem states that, if the metric is of Berwald type (i.e., its Chern--Rund connection defines an affine connection on the underlying manifold), then it is affinely equivalent to a Riemann space, meaning that its affine connection is the Levi--Civita connection of some Riemannian metric. We show for the first time that this result does not extend to general Finsler spacetimes. More precisely, we find a large class of Berwald spacetimes for which the Ricci tensor of the affine connection is not symmetric. The~fundamental difference from positive definite Finsler spaces that makes such an asymmetry possible is the fact that generally, Finsler spacetimes satisfy certain smoothness properties only on a proper conic subset of the slit tangent bundle. Indeed, we prove that when the Finsler Lagrangian is smooth on the entire slit tangent bundle, the Ricci tensor must necessarily be symmetric. For large classes of Finsler spacetimes, however, the Berwald property does not imply that the affine structure is equivalent to the affine structure of a pseudo-Riemannian metric. Instead, the affine structure is that of a metric-affine geometry with vanishing torsion.}


\newcommand{\D}{\text{d}}

\begin{document}

\section{Introduction}\label{sec:intro}
Pseudo-Finsler geometry is a promising mathematical framework for an improved geometric description of the gravitational interaction, beyond the pseudo-Riemannian geometry employed in general relativity~\cite{Bogoslovsky1994,Pfeifer:2019wus,Kostelecky:2011qz,Pfeifer:2011xi,Fuster:2015tua,Gibbons:2007iu,Kouretsis:2008ha,Lammerzahl:2018lhw,Mavromatos:2010nk,Rutz,Stavrinos2014}. Most recently, Finsler spacetime geometry was suggested as the optimal mathematical language to describe the gravitational field of kinetic gases~\cite{Hohmann:2019sni}. While positive definite Finsler geometry, as an extension of Riemannian geometry, is a long-standing, well-established field in mathematics~\cite{Finsler,Bao,Miron}, pseudo-Finsler geometry is still in the process of being developed~\cite{Beem,Lammerzahl:2012kw,Javaloyes:2018lex,Hohmann:2018rpp,Minguzzi:2014fxa,Tavakol2009,TAVAKOL198523}.

The (pseudo-)Finslerian manifolds that can be regarded as closest to the (pseudo-)Riemannian manifolds are the so-called Berwald spaces, resp. spacetimes~\cite{Fuster:2018djw,Szilasi2011,Gallego2017,Berwald1926}. For these geometries, the~canonical nonlinear connection, which is the fundamental building block of the Finsler geometry under consideration, is actually linear in its dependence on the directional variable of the tangent bundle, and thus gives rise to an affine connection on the base manifold. Hence, the natural question to ask is: Under which conditions are Berwald geometries metrizable, i.e., under which conditions does there exist a (pseudo-)Riemannian metric such that the affine Berwald connection is the Levi--Civita connection of this metric?

For Berwald spaces, i.e., in the positive definite case, this question is answered by Szabo's theorem~\cite{Szabo}. It states that every Berwald space is Riemann metrizable. In other words, for every Berwald space, there exists a Riemannian metric such that the affine connection of the Berwald space is identical to the Levi--Civita connection of the metric. The desired Riemannian metric can be constructed explicitly from the Finsler metric of the Berwald space by averaging the Finsler metric over the indicatrix of the Finsler function~\cite{Crampin}. For Finsler spacetimes, no extension of Szabo's theorem has been presented so far. 

In this article, we show for the first time that it is not possible to extend Szabo's theorem to Finsler spacetimes, and thus, Berwald spacetimes are in general non-metrizable. Instead, the following statement holds. The affine connection of every Berwald spacetime is identical to a metric-affine geometry on the base manifold, with a torsion-free connection that is in general not metric compatible. More explicitly, we find examples of Berwald spacetime metrics such that the Ricci tensor of the affine connection is not symmetric. The origin of this lack of symmetry lies in the weaker smoothness assumptions on the Finsler Lagrangian compared to the Finsler function on Finsler spaces.

We present this result as follows. In Section \ref{sec:fins}, we introduce the necessary notions of the geometry of Finsler spaces and Finsler spacetimes before we present our main result, the non-metrizability of Berwald--Finsler spacetimes, in Section \ref{sec:Berwald}. We close this article with a short discussion of our results in Section \ref{sec:disc}.
 
\section{Finsler Geometry}\label{sec:fins}
In Finsler geometry, the geometry of a manifold is derived from a general geometric length measure for curves, defined by a so-called Finsler function. The origin of this idea goes back to Riemann~\cite{Riemann}, but it was only systematically investigated by Finsler in his thesis~\cite{Finsler}. For positive definite geometries, the generalization of Riemannian geometry to Finsler geometry has long been known in the literature. Applying the same kind of generalization from pseudo-Riemannian geometry to pseudo-Finsler geometry is not so straightforward.

We will recall the definition of Finsler spaces and Finsler spacetimes and point out the differences in the construction.

Throughout this article, we consider the tangent bundle $TM$ of an $n$-dimensional manifold $M$, equipped with manifold induced local coordinates, as follows. A point $(x, \dot x)\in TM$ is labeled by the coordinates $(x^a,\dot x^a)$ given by the decomposition of the vector $\dot x = \dot x^a \partial_a \in T_xM$, where $x^a$ are the local coordinates of the point $x \in M$. If there is no risk of confusion, we will sometimes suppress the indices of the coordinates. The local coordinate bases of the tangent and cotangent spaces, $T_{(x,\dot x)}TM$ and $T^*_{(x,\dot x)}TM$, of the tangent bundle are $\{\partial_a = \frac{\partial}{\partial x^a},\dot \partial_a = \frac{\partial}{\partial \dot x^a}\}$ and $\{dx^a, d\dot x^a\}$. In the following, unless elsewhere specified, by smooth, we will always understand $\mathcal{C}^{\infty}$-smooth.

\subsection{Finsler Spaces}\label{ssec:FS}
A Finsler space is a pair $(M,F)$, where $M$ is a smooth $n$-dimensional manifold and the Finsler function $F:TM \to \mathbb{R}$ is continuous on $TM$ and smooth on $TM\setminus\{0\}$, and (see~\cite{Bao}):
\begin{itemize}
	\item $F$ is positively homogeneous of degree one with respect to $\dot x$: $F(x,\alpha \dot x) = \alpha F(x,\dot x)$ for all $\alpha \in \mathbb{R}^+$,
	\item the matrix:
	\begin{align}
	g^F_{ab} = \tfrac{1}{2}\dot \partial_a \dot \partial_b F^2\,
	\end{align}
	defines the Finsler metric tensor $g^F$ and is positive definite at any point of $TM\setminus\{0\}$, in one (and then, in any) local chart around that point.
\end{itemize}

The length of a curve $\gamma:\mathbb{R}\supset [t_1,t_2] \to M$ on a Finsler space is defined by the parametrization invariant length integral:
\begin{align}\label{eq:FinsLength}
S[\gamma] = \int_{t_1}^{t_2}d\tau F(\gamma, \dot \gamma)\,.
\end{align}

To discuss Szabo's theorem later, we need the building blocks of the geometry of a Finsler space, the geodesic spray coefficients:
\begin{align}
G^a = \frac{1}{4}g^{Faq}(\dot x^m \partial_m \dot \partial_q F^2 - \partial_q F^2)\,.
\end{align}

They are defined on each coordinate chart of $TM\setminus \{0\}$ and characterize the geodesic equation, which, for arc length parametrized curves, is:
\begin{align}\label{eq:fgeod}
\ddot x^a + 2 G^a(x,\dot x) = 0\,.
\end{align}

Further objects we need for the arguments in this article are the horizontal derivative operators:
\begin{align}\label{eq:horr}
\delta_a = \partial_a - \dot \partial_a G^b \dot \partial_b\,,
\end{align}
the local coordinate expressions of the Chern--Rund connection coefficients:
\begin{align}\label{eq:CRcoef}
\Gamma^{a}{}_{bc} = \frac{1}{2}g^{F aq}(\delta_b g^F_{cq} + \delta_c g^F_{bq} - \delta_q g^F_{bc} )\,,
\end{align}
the $hh$-Chern--Rund curvature:
\begin{align}\label{eq:CRR}
R^c{}_{adb} = \delta_d \Gamma^{ c}{}_{ab} - \delta_b \Gamma^{ c}{}_{ad} + \Gamma^{ c}{}_{ds} \Gamma^{ s}{}_{ab} - \Gamma^{ c}{}_{bs}\Gamma^{ s}{}_{ad}\,,
\end{align}
and its corresponding horizontal Ricci tensor:
\begin{align}\label{eq:CRRic}
R_{ab} = R^m{}_{amb} = \delta_m \Gamma^{ m}{}_{ab} - \delta_b \Gamma^{ m}{}_{am} + \Gamma^{ m}{}_{ms} \Gamma^{ s}{}_{ab} - \Gamma^{ m}{}_{bs}\Gamma^{ s}{}_{am}\,.
\end{align}

It is important to notice that in general, $R_{ab}$ is not symmetric:
\begin{align}\label{eq:CRRicanti}
R_{ab} - R_{ba} = \delta_a \Gamma^{ m}{}_{bm} - \delta_b \Gamma^{ m}{}_{am} = R^c{}_{dab}\dot x^dC_c\,.
\end{align}

Here, $C_a = C^b{}_{ab}$ are the components of the trace of the Cartan tensor $C_{abc} = \frac{1}{2}\dot{\partial}_ag^F_{bc}$. The last equality in \eqref{eq:CRRicanti} can most easily be proven by introducing the function:
\begin{align}\label{eq:f}
f = \ln \sqrt{|\det g^L_{ab}|}\,.
\end{align}

By direct calculation, one finds,
\begin{align}
	\delta_a f &= \frac{1}{2 |\det g^L_{cd}|} \delta_a |\det g^L_{cd}| = \frac{1}{2}g^{Lmn}\delta_a g^L_{mn} = \Gamma^m{}_{am}, \\
	\dot{\partial}_a f &= \frac{1}{2 |\det g^L_{cd}|} \dot{\partial}_a |\det g^L_{cd}| = \frac{1}{2}g^{Lmn}\dot{\partial}_a g^L_{mn} = C_a\,
\end{align}
and therefore: 
\begin{align}\label{eq:deltaGam}
\delta_a \Gamma^{ m}{}_{bm} - \delta_b \Gamma^{ m}{}_{am} = \delta_a \delta_b f - \delta_b \delta_a f = [\delta_a, \delta_b]f = R^c{}_{dab}\dot x^d\dot{\partial}_cf\,.
\end{align}

An alternative proof of \eqref{eq:deltaGam} can be obtained in terms of the trace of Equation (47) of~\cite{javaloyes2019good} or the use of Equation (3.4.5) of~\cite{Bao}.

The geometric objects introduced in this section make sense on local charts of $TM\setminus\{0\}$. Details~on the construction and properties of these objects in Finsler geometry can be found in~\cite{Bao,Miron}.

Turning to Finsler spacetimes, we will consider the same geometric objects and see that in general, these only make sense on subbundles of $TM\setminus\{0\}$.

\subsection{Finsler Spacetimes}\label{ssec:FST}
To discuss Finsler spacetimes properly, we recall the notion of a conic subbundle $\mathcal{Q}$ of $TM$~\cite{Javaloyes:2018lex}, which is a non-empty open submanifold $\mathcal{Q}\subset TM\backslash \{0\}$, with the following properties:
\begin{itemize}
	\item $\pi(\mathcal{Q})=M$, where $\pi:TM\to M$ is the canonical projection;
	\item {conic property:} if $(x,\dot{x})\in \mathcal{Q}$, then for any $\lambda >0:$ $(x,\lambda \dot{x})\in \mathcal{Q}$.
\end{itemize}

By a Finsler spacetime, we will understand in the following a pair $(M,L)$, where $M$ is a smooth $n$-dimensional manifold and the Finsler Lagrangian $L:\mathcal{A} \to \mathbb{R}$ is a smooth function on a conic subbundle $\mathcal{A}\subset TM$, such that:
\begin{itemize}
	\item $L$ is positively homogeneous of degree two with respect to $\dot x$: $L(x,\lambda \dot x) = \lambda^2 L(x,\dot x)$ for all $\lambda \in \mathbb{R}^+$,
	\item on $\mathcal{A}$, the vertical Hessian of $L$, called the $L$-metric, is nondegenerate,
	\begin{align}\label{g_def}
	g^L_{ab}=\dfrac{1}{2}\dfrac{\partial ^{2}L}{\partial \dot{x}^{a}\partial \dot{x}^{b}}
	\end{align}
	\item there exists a conic subset $\mathcal{T}\subset \mathcal{A}$ such that on $\mathcal{T}$, $L>0$, $g$ has Lorentzian signature $(+,-,-,-)~$ and, on the boundary $\partial \mathcal{T}$, $L$ can be continuously extended as $L|_{\partial \mathcal{T}} = 0$.\footnote{It is possible to formulate this property equivalently with the opposite sign of $L$ and the metric $g^L$ of signature $(-,+,+,+)$. We~fixed the signature and sign of $L$ here to simplify the discussion.}
\end{itemize}

This is a refined version of the definition of Finsler spacetimes in~\cite{Hohmann:2018rpp} and basically covers, if one chooses $\mathcal{A} = \mathcal{T}$, the improper Finsler spacetimes defined in~\cite{Bernal:2020bul}.

The one-homogeneous function $F$, which defines the length measure \eqref{eq:FinsLength}, is derived from the Finsler Lagrange function as $F=\sqrt{|L|}$ and interpreted as the proper time integral of observers. For clarity, we list some important sets appearing on Finsler spacetimes and comment on their meaning:
\begin{itemize}
	\item $\mathcal{A}$: the subbundle where $L$ is smooth and $g^L$ is nondegenerate, with fiber $\mathcal{A}_{x} = \mathcal{A} \cap T_xM$, called~the set of {admissible vectors},
	\item $\mathcal{N}$: the subbundle where $L$ is zero, with fiber $\mathcal{N}_x = \mathcal{N} \cap T_xM$,
	\item $\mathcal{A}_0 = \mathcal{A}\setminus\mathcal{N}$: the subbundle where $L$ can be used for normalization, with fiber $\mathcal{A}_{0x} = \mathcal{A}_0 \cap T_xM$,
	\item $\mathcal{T}$: a maximally connected conic subbundle where $L > 0$ and the $L$-metric exists and has Lorentzian signature $(+,-,-,-)$, with fiber $\mathcal{T}_x = \mathcal{T} \cap T_xM$.
\end{itemize}

A major difference between Finsler spacetimes and Finsler spaces, as defined above, is the existence of these different nontrivial subbundles of $TM\setminus\{0\}$. For Finsler spaces in their classical definition, these bundles become trivial, i.e., $\mathcal{A}=\mathcal{A}_0=TM\setminus\{0\}$, $\mathcal{N} = \{0\}$, and $\mathcal{T}=\emptyset$.

The geometry of Finsler spacetimes is only well defined on $\mathcal{A}$. Some operations, like integration with a canonical zero-homogeneous length measure, can even only be performed on $\mathcal{A}_0$~\cite{Hohmann:2018rpp}. On $\mathcal{A}$, the geodesic spray is given by:
\begin{align}
G^a = \frac{1}{4}g^{Laq}(\dot x^m \partial_m \dot \partial_q L - \partial_q L)\,.
\end{align}

All further geometric objects, which we introduced in the context of Finsler spaces, are defined by the same expressions \eqref{eq:horr}, \eqref{eq:CRcoef}, \eqref{eq:CRR}, \eqref{eq:CRRic}, and \eqref{eq:CRRicanti}, but, just as the geodesic spray, only make sense on $\mathcal{A}$.

The property that in general, the geometry of a Finsler spacetime is not defined on all of $TM\setminus\{0\}$ is crucial in our following finding, that, in general, Szabo's theorem cannot be extended to Finsler~spacetimes.

\section{Berwald Spacetime Geometry and Metric-Affine Spacetime Geometry with Non-Metricity}\label{sec:Berwald}
A Finsler space or Finsler spacetime is said to be of Berwald type, or simply a Berwald space or spacetime, if the geodesic spray is quadratic in $\dot x$:
\begin{align}
G^a(x,\dot x) =\frac{1}{2} G^a{}_{bc}(x)\dot x^b \dot x^c\,.
\end{align}

In this case, the second $\dot x$-derivatives of the geodesic spray coefficients are affine connection coefficients $\dot{\partial}_b\dot{\partial}_cG^a =G^a{}_{bc}(x)$ on the base manifold. 

A standard result in Finsler geometry is that, in general, the geodesic spray can be expressed in terms of the Chern--Rund connection coefficients \eqref{eq:CRcoef} as:
\begin{align}
G^a(x,\dot x) = \frac{1}{2} \Gamma^a{}_{bc}(x,\dot x)\dot x^b \dot x^c\,.
\end{align}

This means that on a Berwald space, or spacetime, the Chern--Rund connection coefficients are independent of $\dot x$ and the affine connection on $M$, defined by the connection coefficients $G^a{}_{bc}(x)$, is precisely given by the Chern--Rund connection, i.e., $\Gamma^a{}_{bc}(x,\dot x) = \Gamma^a{}_{bc}(x) = G^a{}_{bc}(x)$. 

For Berwald spaces, it is known that (see~\cite{Szabo}):
\begin{Theorem}[Szabo's theorem]\label{thm:sz}
	Let $(M,F)$ be a Finsler space of Berwald type. Then, there exists a Riemannian metric $g$ on $M$ such that the affine connection of the Berwald space is the Levi--Civita connection of g
.
\end{Theorem}

Thus, the affine structure of a Berwald space~$(M,F)$ is identical to the affine structure of a Riemannian manifold~$(M,g)$. The metric $g$ can be constructed explicitly from the Finsler metric $g^F$ by an averaging procedure over its $\dot x$ dependence~\cite{Crampin}.

Next we demonstrate that Szabo's theorem can in general not be extended to Berwald spacetimes.

\subsection{A Necessary Condition for the Metrizability of Berwald Spacetimes}
Let $(M,L)$ be a Berwald spacetime with Chern--Rund affine connection coefficients $\Gamma^a{}_{bc} = \Gamma^a{}_{bc}(x)$, as defined in \eqref{eq:CRcoef}. Then, the horizontal Chern--Rund Ricci tensor (see \eqref{eq:CRRic}) is independent of $\dot x$ and takes in every coordinate chart the form:
\begin{align}\label{eq:BerRic}
R_{ab} = \partial_m \Gamma^m{}_{ab} - \partial_b \Gamma^m{}_{am} + \Gamma^m{}_{ms} \Gamma^s{}_{ab} - \Gamma^m{}_{bs} \Gamma^s{}_{am}\,.
\end{align}

It can be regarded as the Ricci tensor of the affine connection with coefficients $\Gamma^a{}_{bc}(x)$ on $M$. A necessary, but not sufficient, condition for the connection defined by $\Gamma^a{}_{bc}$ to be the Levi--Civita connection of a pseudo-Riemannian metric is that the Ricci tensor \eqref{eq:BerRic} is symmetric.

From \eqref{eq:CRRicanti}, we find that for Berwald geometries, the skew-symmetric part is given by:
\begin{align}
R_{ab}(x) - R_{ba}(x) = R^c{}_{dab}(x) \dot x^d C_c(x,\dot x)\,,
\end{align}
where we expressed the explicit dependence on variables $x$ and $\dot x$ to highlight that this equation encodes much information about the geometry of Berwald spacetimes. In particular, it gives rise to the following theorem:

\begin{Theorem}[Symmetric Ricci tensor on smooth Berwald spacetimes]\label{thm:II}
	If $(M,L)$ is a Berwald spacetime with $\mathcal{A} = TM\setminus\{0\}$, then $R_{ba}(x) = R_{ab}(x)$.
\end{Theorem}

The proof of Theorem \ref{thm:II} will be presented in Appendix \ref{app:prfthm1}.

In other words, if $L$ is smooth and $g^L$ is non-degenerate on $TM\setminus\{0\}$, the Ricci tensor of a Berwald--Finsler geometry is symmetric. For Finsler spaces, as we defined them in this article, these conditions are satisfied by definition; see Section \ref{ssec:FS}. Yet, in the class of Finsler spacetimes, there are many interesting classes of examples for which $\mathcal{A}$ is not entirely $TM\setminus\{0\}$, but usually only a subset. This is the origin of the existence of Berwald spacetimes with a non-symmetric Ricci tensor. 

If one weakens the definition of Finsler spaces and allows for Finsler functions that are not smooth everywhere on $TM\setminus\{0\}$, as for example conic Finsler geometries, introduced in~\cite{Javaloyes:2018lex}, with a positive definite Finsler metric, then also these admit examples of Berwald type for which the Ricci tensor is not symmetric and hence provide Berwald spaces that are not metrizable. The same arguments we presented for Finsler spacetimes hold in the positive definite case, since the main point is the non-smoothness of the Finsler function/Finsler Lagrangian on $TM\setminus\{0\}$.

Next, we explicitly present a class of Finsler spacetimes that are, in general, not metrizable, since their Ricci tensor is not symmetric.

\subsection{Non-Metrizable Berwald--Finsler Spacetimes}
The following Finsler spacetimes are of Berwald type, but not metrizable, since their Ricci tensor is not symmetric. As we have seen in the previous section, the reason is that the set of admissible vectors $\mathcal{A}$ is smaller than $TM\setminus\{0\}$, which is rather the usual case for Finsler spacetimes, than the exception.

Let $\alpha(\dot x,\dot x) = \alpha_{ab}(x)\dot x^a \dot x^b$ be constructed from a pseudo-Riemannian metric $\alpha$, $\beta(\dot x) = \beta_a(x)\dot x^a$ a one-form on $M$ evaluated on a generic tangent vector, and $c,m,p$ real numbers. Moreover, define the zero-homogeneous variable $s= \frac{\beta(\dot x)^2}{\alpha(\dot x, \dot x)}$. Consider the $(\alpha,\beta)$-Finsler Lagrangian:
\begin{align}\label{eq:Lagr_ext_VGR}
L(x,\dot x) = \alpha(\dot x, \dot x) s^{-p} (c +m\ s )^{p+1}\,.
\end{align}
A brief discussion about the causal properties of this type of Finsler Lagrangians can be found in  Appendix \ref{app:BogKrop}.

Finsler spacetimes $(M,L)$ are generalizations of Bogoslovsky/p-Kropina/very general relativity geometries~\cite{Bogoslovsky1994,Kouretsis:2008ha,Fuster:2015tua,Gibbons:2007iu,Fuster:2018djw,Cohen:2006ky,Kropina,Gomez-Lobo:2016qik,Elbistan:2020mca}, which we recover for $c=1$, $m=0$. The classical Kropina case is included by setting $p=1$. In \cite[Corollary 3]{Pfeifer:2019tyy}
, we have shown that $(M,L)$ is a Berwald spacetime if and only if the covariant derivative of the one-form $\beta$ with respect to the Levi--Civita connection of $\alpha$ satisfies:
\begin{align}\label{eq:Berwald_condition}
\nabla_a\beta_b = H \left([c (1 - p ) + m \alpha^{-1}(\beta,\beta)]\beta_a \beta_b + c p \alpha^{-1}(\beta,\beta)\alpha_{ab}\right)\,,
\end{align}
for an arbitrary function $H= H(x)$ on $M$. The resulting geodesic spray is:
\begin{align}
G^a(x,\dot x) = \frac{1}{2} \Gamma^a{}_{bc}(x) \dot x^b \dot x^c = \frac{1}{2} \big( \gamma^a{}_{bc}(x) -H \left( p c ( \delta^a_{b} \beta_{c} + \delta^a_{c} \beta_{b} ) - \beta^a ( m \beta_b \beta_c + p c\alpha_{bc} ) \right) \big) \dot x^b \dot x^c\,,
\end{align}
and thus, the resulting affine connection coefficients are:
\begin{align}
\Gamma^a{}_{bc}(x) = \gamma^a{}_{bc}(x) -H \left( c p ( \delta^a_{b} \beta_{c} + \delta^a_{c} \beta_{b} ) - \beta^a ( m \beta_b \beta_c + c p \alpha_{bc}) \right)\,.
\end{align}

Here, $\gamma^a{}_{bc}$ are the Christoffel symbols of the pseudo-Riemannian metric $\alpha$; $R[\gamma]_{ab}$ below is the corresponding Ricci tensor. For the Chern--Rund Ricci tensor \eqref{eq:BerRic}, we find:

\begin{align}
R_{bd} 
&= R[\gamma]_{bd} + c p \alpha_{bd} \big( H^2 \alpha^{-1}(\beta,\beta)( c + 3cp + m \alpha^{-1}(\beta,\beta)) + \beta^a \partial_aH\big) \nonumber\\
&+ \beta_b \beta_d \big(2 c p H^2 (c + m \alpha^{-1}(\beta,\beta)) + m \beta^a \partial_aH \big) \nonumber\\
&- \beta_b \partial_dH (m \alpha^{-1}(\beta,\beta)-3 c p ) - c p \beta_d\partial_bH\,.
\end{align}

Thus, for the skew-symmetric part,
\begin{align}\label{eq:Antisymmetric_part_Ricci}
\frac{1}{2}\left(R_{ab} - R_{ba}\right) = \frac{1}{2}( 4c p - m \alpha^{-1}(\beta,\beta) )(\beta_a \partial_bH - \beta_b \partial_aH)\,.
\end{align}

We summarize these findings as follows.

\begin{Proposition}
	Let $(M,L)$ be a Finsler spacetime defined by a Finsler Lagrangian of the type $L(x,\dot x) = \alpha(\dot x, \dot x) s^{-p} (c +m\ s )^{p+1}$ as in Equation \eqref{eq:Lagr_ext_VGR}. Then, the skew-symmetric part $\mathcal{A}(R)$ of the Chern--Rund Ricci tensor is given by:
	\begin{align}
		\mathcal{A}(R) = f \beta \wedge dH\,,
	\end{align}
	where $f=\frac{1}{2}(4cp-m\alpha^{-1}(\beta,\beta))$ is a scalar function. If $f\neq0$ and $\beta \wedge dH \neq0$, then $(M,L)$ is non-metrizable. 
\end{Proposition}

This shows that the Ricci tensor corresponding to the Finsler Lagrangian \eqref{eq:Lagr_ext_VGR} may not be symmetric when $H \neq 0$, in particular when the one-form $\beta$ is not covariantly constant (cf.\eqref{eq:Berwald_condition}). This is indeed often the case\footnote{We will elaborate on this in forthcoming work.}, for example in the class of $(\alpha,\beta)$-Kundt spacetimes introduced in~\cite{Pfeifer:2019tyy}, with Lagrangian \eqref{eq:Lagr_ext_VGR}. 
For the sake of conciseness, we present below a very simple geometry in this class for which the Ricci tensor is {not} symmetric, providing an explicit counterexample that shows that Szabo's theorem does not hold in general for Berwald spacetimes.
	
Let the Lorentzian metric $\alpha$ and the one-form $\beta$ be given by:
\begin{align}
\alpha = 2\, \D u \D v + v \,\phi(x,y)\,\D u^2 + \D x^2 + \D y^2,\qquad \beta = \D u,
\end{align}
where $u$, $v$ are so-called light-cone coordinates and $\phi$ is a scalar function. It is easy to see that $\beta$ is a null one-form with respect to $\alpha$, i.e., $\alpha^{-1}(\beta,\beta)=0$. Assuming $c\neq 0$ and $p\neq 1$, the pair $\alpha,\beta$ satisfies condition \eqref{eq:Berwald_condition} with $H = \frac{\phi}{2c(p-1)}$; hence, the resulting Finsler Lagrangian $L = \alpha s^{-p} (c +m s )^{p+1}$ is of Berwald type. 

Plugging the function $H$ in Equation \eqref{eq:Antisymmetric_part_Ricci}, we find the following non-trivial components of the skew-symmetric part of the Ricci tensor:
\begin{align}
\frac{1}{2}\left(R_{ux}-R_{xu}\right) = \frac{p}{p-1}\partial_x\phi, \\
\frac{1}{2}\left(R_{uy}-R_{yu}\right) = \frac{p}{p-1}\partial_y\phi.
\end{align}

We see that they are independent of the parameters $c$ and $m$ and do not vanish for non-constant $\phi$ and $p\neq 0$. This simple counterexample shows that the affine connection of a Berwald spacetime $(M,L)$ is not necessarily equivalent to the Levi--Civita connection of a (pseudo-)Riemannian metric. Hence, Szabo's theorem does not extend to general Finsler spacetimes.

We would like to point out that if one would exchange the Lorentzian metric in the examples with a Riemannian metric, the same derivations can be made. Hence, if one allows for Finsler spaces with a Finsler function that is not smooth on all of $TM\setminus\{0\}$, as for instance conic Finsler geometries~\cite{Javaloyes:2018lex}, one~also finds counterexamples to Szabo's theorem for positive definite Finsler geometries.

\subsection{Affine Structure of Berwald Spacetimes}\label{sec:affine}
In physics, metric affine theories of gravity are considered as an extension of general \mbox{relativity~\cite{Hehl:1994ue,Percacci:2019hxn}} and an effective description of quantum gravity~\cite{Percacci:2020bzf}. In these theories, a spacetime manifold is equipped with a spacetime metric and an independent affine connection, which in general is not metric compatible, possesses curvature and torsion. Special instances of metric affine gravity are those where the connection has either only torsion~\cite{Krssak:2018ywd}, is only not-metric compatible~\cite{Lu:2019hra,Nester:1998mp,Barros:2020bgg}, has only curvature (the usual pseudo-Riemannian case), or possesses any possible combination of these three properties. In this section we briefly summarize that the geometry of non-metrizable Berwald spacetimes is equivalent to a metric affine geometry with a torsion free, non-metric compatible~connection.

Pick any pseudo-Riemannian metric $g$, and let $\gamma^a{}_{bc}[g]$ be the Christoffel symbols of its Levi--Civita connection. Then:
\begin{align}
\Gamma^a{}_{bc}(x) = \gamma^a{}_{bc}[g] + D^a{}_{bc}\,,
\end{align}
where $D$ is a $(1,2)$-tensor field on $M$. By construction, the $\Gamma^a{}_{bc}(x)$ are symmetric in their lower indices, and hence, they define a torsion-free affine connection. Thus, by the decomposition of affine connections into Levi--Civita, contorsion, and non-metricity parts, the tensor $D^a{}_{bc}$ defines the non-metricity $Q_{abc} = \nabla_a g_{bc} = - D^s{}_{ac}g_{sb} - D^s{}_{ab}g_{sc}$ of the connection. In view of this argument, we can formulate the following proposition:

\begin{Proposition}
	Let $(M,L)$ be a Berwald--Finsler spacetime, and let $\Gamma^a{}_{bc}$ be the induced affine connection coefficients on~$M$. Moreover, choose any pseudo-Riemannian metric $g$. The affine structure of the Finsler spacetime $(M,L)$ is equivalent to the affine structure of the metric-affine geometry of $(M,g,\Gamma)$, where the connection defined by the connection coefficients $\Gamma^a{}_{bc}$ is torsion free, but in general not metric compatible. A~Berwald--Finsler spacetime is metrizable if and only if there exists a pseudo-Riemannian metric $g$ such that the non-metricity $Q_{abc}$ vanishes.
\end{Proposition}

\section{Discussion}\label{sec:disc}
Evidently, Finsler spacetime geometry has a very different behavior compared to that of Finsler spaces. The origin of this difference lies, on the one hand, in the weaker smoothness assumptions on the defining Finsler Lagrangian $L$ and, on the other hand, in the fact that its indicatrix is necessarily non-compact in the indefinite case~\cite{Voicu_2017}. As a consequence, classical theorems that hold on Finsler spaces may not hold on Finsler spacetimes anymore. It was already known that Deicke's theorem~\cite{Deicke1953} does not hold on generic Finsler spacetimes~\cite{Voicu_2017}. With this article, we demonstrated for the first time that the same is true for Szabo's theorem. 

These findings call for a systematic study and classification of the geometric properties of Finsler spacetimes in general and of Berwald spacetimes in particular.

\authorcontributions{The authors all contributed substantially to the derivation of the presented results, as well as the analysis, drafting, review, and finalization of the manuscript. All authors read and agreed to the published version of the manuscript.}

\funding{C.P. was supported by the Estonian Ministry for Education and Science through the Personal Research Funding Grant PSG489, as well as the European Regional Development Fund through the Center of Excellence TK133
 ``The Dark Side of the Universe''. The work of A.F. is part of the research program of the Foundation for Fundamental Research on Matter (FOM
), which is financially supported by the Netherlands Organisation for Scientific Research (NWO).}

\acknowledgments{The authors would like to acknowledge networking support by the COST Action QGMM
 (CA18108), supported by COST (European Cooperation in Science and Technology).}

\conflictsofinterest{The authors declare no conflict of interest.}

\appendixtitles{yes} 
\appendix

\section{Proof of Theorem 2}\label{app:prfthm1}
To identify the origin of the lack of symmetry of the Chern--Rund--Ricci tensor, which implies non-metrizability, we used Theorem~\ref{thm:II}, which we state here again:

\begin{center}
	\textit{If $(M,L)$ is a Berwald spacetime with $\mathcal{A} = TM\setminus\{0\}$, then $R_{ba}(x) = R_{ab}(x)$.}
\end{center}

We now provide the proof of the theorem. The starting point of the proof is Equation \eqref{eq:CRRicanti}:
\begin{align}
R_{ab} - R_{ba} = \delta_a \Gamma^{ m}{}_{bm} - \delta_b \Gamma^{ m}{}_{am} = R^c{}_{dab}\dot x^dC_c\,,
\end{align}
which on Berwald spacetimes implies that the expression:
\begin{align}
\phi(x,\dot x) = R^c{}_{dab}(x)\dot x^dC_c(x,\dot x) = R^c{}_{dab}(x)\dot x^d\dot{\partial}_cf(x,\dot x)
\end{align}
is actually independent of $\dot x$. It is clear that the Ricci tensor of a Berwald spacetime is symmetric if and only if $\phi$ identically vanishes, which we will now connect to the zeros of the derivatives of the function $f$. The function $f = \ln \sqrt{|\det g^L_{ab}|}$, as defined in \eqref{eq:f}, is zero-homogeneous in its dependence on $\dot x$, which~means it naturally lives on the positive projective tangent bundle $PTM^+$; see~\cite{Hohmann:2018rpp}. 

Recall that $PTM^+$ is defined as the set of equivalence classes $[(x,\dot x)]$, where $(x, \dot x) \sim (x', \dot x')$ if and only if $ (x', \dot x') = (x,\lambda \dot x)$ for some positive real $\lambda$. This makes $PTM^+$ a $2n-1$ dimensional manifold with coordinate charts built as follows. Consider a coordinate chart $(U,\varphi)$ on $M$ and define the open subsets $V^+_i = \{(x,\dot x)\in TU | \dot x^i > 0\}$ and $V^-_i = \{(x,\dot x)\in TU | \dot x^i < 0\}$ on $TM$. Then, for each $[(x,\dot x)]$ with $(x,\dot x)\in V^\pm_i$, we define the coordinates: 
\begin{align}
(x^a, u^\alpha) = \left(x^0, ... x^n, \frac{\dot x^0}{\dot x^i}, ..., \frac{\dot x^{i-1}}{\dot x^i}, \frac{\dot x^{i+1}}{\dot x^i},...,\frac{\dot x^{n-1}}{\dot x^i}\right)\,,
\end{align}
where $a=0,...,n-1$ and $\alpha = 0,...,n-2$. More conveniently, one can use so-called homogeneous coordinates $(x^a,\dot x^a)$, which are nothing but the coordinates on $TM$ of an arbitrary representative of the equivalence class $[(x,\dot x)]$. Homogeneous coordinates are only defined up to a scaling factor. $PTM^+$~itself is a fiber bundle over $M$, with compact fibers, diffeomorphic to Euclidean spheres~\cite{Chern-Shen-Lam}. Functions on $TM\setminus\{0\}$ can be understood as functions on $PTM^+$ if and only if they are zero-homogeneous in $\dot x$; in homogeneous coordinates, calculus on $PTM^+$ is formally identical to the one on $TM$.

Assume that $\mathcal{A} = TM\setminus \{0\}$. This implies that $|\det g|$ is smooth and nonzero on $TM\setminus \{0\}$, and so is $f$. Moreover, since they are zero-homogeneous in $\dot x$, they can be regarded as functions on $PTM^+$.

Fix an arbitrary $x\in M$ and an arbitrary local chart around $x$, then $f_x(\cdot) = f(x,\cdot)$ is defined on the fiber $PT_xM^+$, which is compact. Since $f_x$ is smooth, it admits at least a local extremum, say at $\dot x = v$, and hence, $\dot{\partial}_a f_x|_{\dot x =v} = 0$. For the function $\phi(x,\dot x)$, this implies:
\begin{align}
\phi(x,v) = R^c{}_{dab}(x)v^d\dot{\partial}_cf(x,v) = 0\,.
\end{align}

On the other hand, on Berwald spacetimes, $\phi$ is independent of $\dot x$, and we can conclude that $\phi$ is identically zero for all $\dot x$, which completes the proof.

As a remark, the crucial ingredient of the proof is that $|\det g^L_{ab}|$ is smooth and nonzero on the entire $TM\setminus \{0\}$. If $|\det g^L_{ab}|$ were only smooth and nonzero on a smaller subset, then $f_x$ would not be defined on the entire fiber $PT_xM^+$, and the conclusion of the theorem would fail.

\section{Generalized Bogoslovsky/Kropina--Finsler Lagrangians}\label{app:BogKrop}
The class of Finsler Lagrangians we considered as a counterexample to Szabo's theorem are given in Equation \eqref{eq:Lagr_ext_VGR}, which we rewrite here in the best way to briefly investigate the causal structure of the Finsler Lagrangian:
\begin{align}\label{app:Lagr_ext_VGR}
L(x,\dot x) = \alpha(\dot x, \dot x) s^{-p} (c +m\ s )^{p+1} 
			&= \frac{\alpha(\dot x,\dot x)^{p+1}}{\beta(\dot x)^{2p}}\left(c + m \frac{\beta(\dot x^2)}{\alpha(\dot x, \dot x)}\right)^{p+1}\\
			&= \frac{\zeta(\dot x, \dot x)^{p+1}}{\beta(\dot x)^{2p}}\,,
\end{align}
where we introduced the shorthand notation of an effective bilinear form $\zeta(\dot x, \dot x) = c \alpha(\dot x, \dot x) + m \beta(\dot x)^2$. Hence, the Finsler Lagrangian in consideration is effectively of Kropina/Bogoslovsky/VGR
, for which the metric $\zeta$ is constructed from more fundamental building blocks. Depending on the properties of the building blocks $\beta$, $\alpha$, $m$, and $c$, $\zeta$ can have different signatures.

The causal structure of $L$ can be characterized for three different ranges of values for the parameter~$p$:
\begin{enumerate}
	\item $p>0$: $L=0\Leftrightarrow \zeta(\dot x, \dot x) = 0$, and $L$ is not defined for $\beta(\dot x) = 0$;
	\item $0>p>-1$: $L=0\Leftrightarrow \zeta(\dot x, \dot x) = 0$ or $\beta(\dot x) = 0$;
	\item $p<-1$: $L=0\Leftrightarrow \beta(\dot x) = 0$, and $L$ is not defined for $\zeta(\dot x, \dot x) = 0$.
\end{enumerate}

The third case never leads to a Finsler spacetime since the null set $\beta(\dot x)=0$ singles out a hyperplane and never allows for the existence of a convex cone of timelike vectors.

For the other two cases, a necessary condition to obtain a Finsler spacetime is that the bilinear form $\zeta$ is a pseudo-Riemannian metric of a Lorentzian signature \cite[Appendix B]{Hohmann:2018rpp}, where the same class of Finsler Lagrangians is studied for $p=-q$. This demand leads to conditions on $\beta$, $\alpha$, $m$, and $c$, from~the determinant: 
\begin{align}
\det \zeta_{ab} = c^3 \det (\alpha_{ab})(c + m \alpha^{-1}(\beta,\beta))\,,
\end{align}
which must be negative.

For a choice of $\beta$, $\alpha$, $m$, and $c$ such that this necessary requirement is satisfied, one can apply the classification done in \cite[Appendix B]{Hohmann:2018rpp} to identify viable Finsler spacetimes. A main finding there is that for $-1 < p < 1$ and $\beta$ being $\zeta$ timelike, the cone of future pointing timelike vectors $\mathcal{T}$ is given by the cone of future pointing timelike vectors of $\zeta$.

A complete classification of the Finsler Lagrangians \eqref{eq:Lagr_ext_VGR} goes beyond the scope of this article and is left for future investigation.

\reftitle{References}



\end{document}